\documentclass{gtart_a}
\pdfoutput=1

\title{Thin buildings}

\author{Jan Dymara}
\givenname{Jan}
\surname{Dymara}
\address{Instytut Matematyczny\\
Uniwersytet Wroc\l awski\\\newline 
pl. Grunwaldzki 2/4\\
50-384 Wroc\l aw\\Poland}
\email{dymara@math.uni.wroc.pl}
\urladdr{}

\volumenumber{10}
\issuenumber{}
\publicationyear{2006}
\papernumber{18}
\lognumber{0698}
\startpage{667}
\endpage{694}

\doi{}
\MR{}
\Zbl{}


\received{6 January 2006}
\revised{}
\accepted{30 April 2006}
\published{24 May 2006}
\publishedonline{24 May 2006}
\proposed{Wolfgang L\"uck}
\seconded{Martin Bridson, Steve Ferry}
\corresponding{}
\editor{}
\version{}

\AtBeginDocument{\let\bar\wbar}
\numberwithin{equation}{section}
\allowdisplaybreaks
\def\Aut{\mathrm{Aut}}

\makeatletter
\def\cnewtheorem#1[#2]#3{\newtheorem{#1}{#3}[section]
\expandafter\let\csname c@#1\endcsname\c@thm}

\theoremstyle{plain}
\newtheorem{thm}{Theorem}[section]
\cnewtheorem{propo}[thm]{Proposition}
\cnewtheorem{coro}[thm]{Corollary}
\cnewtheorem{lem}[thm]{Lemma}
\newtheorem*{claim}{Claim}
\theoremstyle{definition}
\newtheorem*{definition}{Definition}
\theoremstyle{remark}
\newtheorem*{remark}{Remark}

\makeatother  

\newcommand{\supp}{{\rm supp}}
\newcommand{\metric}{Moussong metric }
\newcommand{\rhtm}{right Hecke $t$--mul\-ti\-pli\-ca\-tion }
\newcommand{\ce}[1]{\langle #1\rangle }

\begin{document}

\begin{asciiabstract}
Let X be a building of uniform thickness q+1.  L^2-Betti numbers of X
are reinterpreted as von-Neumann dimensions of weighted L^2-cohomology
of the underlying Coxeter group.  The dimension is measured with the
help of the Hecke algebra.  The weight depends on the thickness q.
The weighted cohomology makes sense for all real positive values of q,
and is computed for small q.  If the Davis complex of the Coxeter
group is a manifold, a version of Poincare duality allows to deduce
that the L^2-cohomology of a building with large thickness is
concentrated in the top dimension.
\end{asciiabstract}

\begin{htmlabstract}
Let X be a building of uniform thickness q+1.  L<sup>2</sup>&ndash;Betti
numbers of X are reinterpreted as von-Neumann dimensions of weighted
L<sup>2</sup>&ndash;cohomology of the underlying Coxeter group. The
dimension is measured with the help of the Hecke algebra. The weight
depends on the thickness q. The weighted cohomology makes sense for all
real positive values of q, and is computed for small q. If the Davis
complex of the Coxeter group is a manifold, a version of Poincar&eacute;
duality allows to deduce that the L<sup>2</sup>&ndash;cohomology of
a building with large thickness is concentrated in the top dimension.
\end{htmlabstract}

\begin{abstract}
Let $X$ be a building of uniform thickness $q+1$.  $L^2$--Betti numbers
of $X$ are reinterpreted as von-Neumann dimensions of weighted
$L^2$--cohomology of the underlying Coxeter group. The dimension is
measured with the help of the Hecke algebra. The weight depends on the
thickness $q$. The weighted cohomology makes sense for all real
positive values of $q$, and is computed for small $q$. If the Davis
complex of the Coxeter group is a manifold, a version of Poincar\'e
duality allows to deduce that the $L^2$--cohomology of a building with
large thickness is concentrated in the top dimension.
\end{abstract}

\maketitle

\section*{Introduction}

Let $(G,B,N,S)$ be a $BN$--pair, and let $X$ be the associated building
(notation as in Brown \cite[Chapter 5]{Br}).  There are many geometric
realizations of $X$. We consider the one introduced by Davis in
\cite{D1}. Then $X$ is a locally finite simplicial complex, acted upon
by $G$. The action has a fundamental domain with stabiliser $B$.  The
standard choice of such a domain is called the Davis chamber.  We can
and will assume that $G$ is a closed subgroup of the group $\Aut(X)$ of
simplicial automorphisms of $X$ (in the compact-open topology). If
this is not the case, one can pass to the quotient of $G$ by the
kernel of the $G$--action on $X$ (that quotient is a subgroup of
$\Aut(X)$), and then take its closure in $\Aut(X)$.
 
Let $L^2C^i(X)$ be the space of $i$--cochains on $X$ which are
square-summable with respect to the counting measure on the set
$X^{(i)}$ of $i$--simplices in $X$.  Then the coboundary map
$\delta^i\colon L^2C^i(X)\to L^2C^{i+1}(X)$ is a bounded operator. The
reduced $L^2$--cohomology of $X$ is defined to be
$L^2H^i(X)=\ker{\delta^i}/\overline{{\rm im}{\delta^{i-1}}}$.  This is
a Hilbert space, carrying a unitary $G$--representation.  Using the von
Neumann $G$--dimension one defines $L^2b^i(X)=\dim_GL^2H^i(X)$.  We are
interested in calculating these Betti numbers.  (This problem was
considered by Dymara and Januszkiewicz in \cite{DJ} and by Davis and
Okun in \cite{DO}.)

The first step is to pass from the cochain complex
$(L^2C^*(X),\delta)$ to a smaller complex of $B$--invariants:
$(L^2C^*(X)^B,\delta)$. Now $L^2C^i(X)^B$ can be 
identified with a space of cochains on $X/B=\Sigma$---the Davis complex
of the Weyl group $W$ of the building.
However, a simplex  $\sigma\in\Sigma$ has a preimage in $X$ consisting
of $q^{d(\sigma)}$ simplices, where
$q+1$ is the thickness of the building and $d(\sigma)$
is the distance from $\sigma$ to the chamber stabilised by $B$.
Therefore a cochain $f$ on $\Sigma$ represents 
a square-summable $B$--invariant cochain
if and only if it satisfies
$\sum_\sigma |f(\sigma)|^2q^{d(\sigma)}<\infty$;
we denote the space of such cochains $L^2_qC^*(\Sigma)$.
The complex $(L^2_qC^*(\Sigma),\delta)$ 
and its (reduced) cohomology $L^2_qH^*(\Sigma)$ 
are acted upon by the Hecke
algebra ${\bf C}[B\backslash  G/B]$. A suitable von Neumann completion 
of the latter can be used to measure the dimension 
of $L^2_qH^i(\Sigma)$, yielding Betti numbers
$L^2_qb^i(\Sigma)$. It turns out that        
$L^2_qb^i(\Sigma)=L^2b^i(X)$. In particular, 
the $L^2$--Betti numbers of a building depend
only on $W$ and on $q$.

The good news is that the complex
$(L^2_q(\Sigma),\delta)$, the Hecke algebra and the Betti
numbers $L^2_qb^i(\Sigma)$ can be defined for all real $q>0$,
in a uniform manner which for integer values of $q$
gives exactly the objects discussed above.
It turns out that for small $q$ (namely
for $q<\rho_W$, where $\rho_W$ is the logarithmic growth
rate of $W$) the Betti numbers $L^2_qb^i(\Sigma)$ 
are $0$ except for $i=0$. Since $\rho_W\le1$, this result
says nothing about actual buildings. 
However, in \fullref{sec:6} we prove a version of Poincar\'e
duality, saying that if $\Sigma$ is a manifold of dimension $n$, 
then $L^2_qb^i(\Sigma)=L^2_{1/q}b^{n-i}(\Sigma)$. 
Thus, if the Davis complex of the Weyl group of a building 
(ie, an apartment in the Davis realization
of the building)
is an $n$--manifold, and if $q>{1\over\rho_W}$, then 
the $L^2$--Betti numbers of the building vanish except for 
$L^2b^n(X)$.

Examples of buildings to which our method  applies can be constructed
from flag triangulations of spheres. Davis associates a right-angled
Coxeter group to any such triangulation; this right-angled Coxeter group 
is the Weyl group of a family 
of buildings with manifold apartments, parametrised by thickness.
Let us mention that the argument applies also to Euclidean buildings,
yielding another calculation of their 
$L^2$--Betti numbers. 

In a forthcoming paper (Davis--Dymara--Januszkiewicz--Okun
\cite{DDJO}) the $L^2$--Betti numbers of all buildings satisfying
$q>{1\over\rho_W}$ are calculated.

The definitions, results and arguments of this paper
go through, with appropriate reading, in the multi-parameter case.
A detailed account of the multi-parameter setting is given in \cite{DDJO}.

The author thanks Dan Boros, 
Tadeusz Januszkiewicz, Boris Okun
and especially Mike Davis 
for useful discussions.

The author was partially supported by KBN grants 5 P03A 035 20
and 2 P03A 017 25,
and by a scholarship of the Foundation for Polish Science.

\setcounter{section}{-1}

\section{Integer thickness}\label{sec:0}

Let $(W,S)$ be a Coxeter system. Let $\Delta$ be a simplex with 
codimension 1 faces labelled by elements of $S$, 
and let $\Delta'$ be its first barycentric subdivision.
Each $T\subseteq S$ generates a subgroup $W_T$ of $W$ called a special 
subgroup; also, $T$ corresponds to a face $\Delta_T$ of $\Delta$
(the intersection of codimension 1 faces labelled by elements of $T$).
The Davis chamber $D$ is the subcomplex of $\Delta'$ spanned by 
barycentres of faces $\Delta_T$ for which $W_T$ is finite
($\cal F$ will denote the set of subsets $T\subseteq S$ such that $W_T$ is finite). 
To every $T\subseteq S$
we assign a face of the Davis chamber: $D_T=D\cap\Delta_T$.
The Davis realization  $\Sigma$ of the Coxeter complex is $W\times D/\sim$,
where $(w,p)\sim(u,q)$ if and only if for some $T$ we have $p=q\in D_T$ 
and $w^{-1}u\in W_T$. 
The action of $W$ on the first factor descends
to an action on $\Sigma$.
We denote the image of
$\sigma$ under the action of $w$ by $w\sigma$, 
and the $W\!$--orbit of $\sigma$ in $\Sigma$ by $W\sigma$. 
The images of $w\times D$  in $\Sigma$ 
are called chambers. The action of $W$ on $\Sigma$ is simply 
transitive on the set of chambers. 

A Tits building $X_{Tits}$ with Weyl group $W$
is a set with a $W\!$--valued distance function
$d$, satisfying certain conditions (see Ronan \cite{Ron}). Its Davis
incarnation is $X=X_{Tits}\times D/\sim$, where
$(x,p)\sim(y,q)$ if and only if for some $T$ we have $p=q\in D_T$ 
and $d(x,y)\in W_T$. 
The images of $x\times D$ in $X$
are called chambers.
         
We will consider only buildings of uniformly bounded thickness,
ie, such that for some constant $N>0$, any $s\in S$ and any $x\in X_{Tits}$
there are no more than $N$ elements $y\in X_{Tits}$ satisfying $d(x,y)=s$.
If this number of $s$--neighbours of $x$ is equal to $q$
for all pairs $(x,s)$, then we say that the building has uniform thickness 
$q+1$. We denote such building $X(q)$ 
(for a right-angled Coxeter group it is unique).

Uniformly bounded thickness is equivalent to $X$ being 
uniformly locally finite. Thus we can consider 
(reduced) $L^2$--(co)homology of $X$. This is obtained from
the complex of $L^2$ (co)chains on $X$ with the usual
(co)boundary operators $\partial$, $\delta$. These operators are
in fact adjoint to each other, so that
the (co)homology can be identified with $L^2{\cal H}^*(X)$, the
space of harmonic (co)chains (``reduced'' means 
that we divide the kernel by the closure of the image).

Assume now that $X_{Tits}$ comes from a $BN$--pair in a group $G$.
Then $G$ acts by simplicial automorphisms on $X$.
We can assume that $G$ acts faithfully and is locally compact 
(possibly taking the closure of its image
in $\Aut(X)$ in the compact-open topology). We use $G$
to measure the size of $L^2{\cal H}^i(X)$ via the von Neumann 
dimension. To do this, we first express 
$L^2C^i(X)$ as $\oplus_{\sigma^i\subset D}L^2(G\sigma^i)$. Then we notice that
$L^2(G\sigma^i)$ is naturally isomorphic to $L^2(G)^{G_{\sigma^i}}$
(where $G_{\sigma^i}$ is the stabiliser of $\sigma^i$ in $G$).
It is convenient to multiply this isomorphism by a suitable scalar 
factor in order to make it isometric.
Then the space $L^2(G)^{G_{\sigma^i}}$ is embedded into
$L^2(G)$, giving us finally an embedding of left $G$--modules
$L^2C^i(X)\hookrightarrow\oplus_{\sigma^i\subset D}L^2(G)$.
In particular,
$L^2{\cal H}^i(X)$ is now embedded as a left $G$--module
in $\oplus_{\sigma^i\subset D}L^2(G)$; we can consider the orthogonal
projection onto this subspace, and define $L^2b^i(X)$ to be the  
von Neumann trace of that projection. 
Let $B$ be the stabiliser of $D$ in $G$. 
For each $\sigma^i\subset D$ we have a vector ${\bf 1}_\sigma$ in
$\oplus_{\sigma^i\subset D}L^2(G)$, having $\sigma$th component
${\bf 1}_B$ and other components 0.
The projection onto $L^2{\cal H}^i(X)$ is given by a matrix whose 
$\sigma$th row gives the projection of
${\bf 1}_\sigma$ on $L^2{\cal H}^i(X)$,
expressed as an element of $\oplus_{\sigma^i\subset D}L^2(G)$
(while applying this matrix we understand multiplication as convolution).
Notice that both ${\bf 1}_\sigma$ and the space  
$L^2{\cal H}^i(X)$ are $B$--invariant; so therefore will be 
the projection of   ${\bf 1}_\sigma$ on $L^2{\cal H}^i(X)$.

\section{Real thickness}\label{sec:1}

For a $w\in W$ we denote by $d(w)$ the length of a shortest word in the 
generators $S$ representing $w$. For a chamber $c=w\times D$ of $\Sigma$
we put $d(c)=d(w)$. For every simplex $\sigma\subset\Sigma$
there is a unique chamber $c\supseteq\sigma$ with smallest $d(c)$;
we put $d(\sigma)=d(c)$. 

For a real number $t>0$ we equip the set $\Sigma^{(i)}$ of $i$--simplices in
$\Sigma$ with the measure
$\mu_t(\sigma)=t^{d(\sigma)}$. We also pick (arbitrarily) orientations of 
simplices in $D$, and extend them $W\!$--equivariantly
to orientations of all simplices in $\Sigma$. This allows us to identify
chains and cochains with functions. We put
$$L^2_tC^i(\Sigma)=L^2_tC_i(\Sigma)=L^2(\Sigma^{(i)}, \mu_t).$$
We now define $\delta^i\colon L^2_tC^i(\Sigma)\to L^2_tC^{i+1}(\Sigma)$ by
$$\delta^i(f)(\tau^{i+1})=\sum_{\sigma^i\subset\tau^{i+1}}
[\tau:\sigma]f(\sigma)$$
and  $\partial_i^t\colon L^2_tC_i(\Sigma)\to L^2_tC_{i-1}(\Sigma)$ by
$$\partial_i^t(f)(\eta^{i-1})=\sum_{\sigma^i\supset\eta^{i-1}}
[\eta:\sigma]t^{d(\sigma)-d(\eta)}f(\sigma)$$
(here $[\alpha:\beta]=\pm1$ tells us whether 
orientations of $\alpha$ and $\beta$ agree or not).
We have 
{\setlength\arraycolsep{2pt}
\begin{eqnarray}
\langle \delta^i(f),g\rangle_t & = & \sum_{\tau^{i+1}}\Bigl(
\sum_{\sigma^i\subset\tau^{i+1}}[\tau:\sigma]f(\sigma)\overline{g(\tau)}
t^{d(\tau)}\Bigr)\nonumber\\
& = &
\sum_{\sigma^i}f(\sigma)\overline{\Bigl(
\sum_{\tau^{i+1}\supset\sigma^i}[\tau:\sigma]t^{d(\tau)-d(\sigma)}g(\tau)
\Bigr)}
t^{d(\sigma)}
=\langle f,\partial_i^t(g)\rangle_t.\nonumber
\end{eqnarray}
}%
That is, $\delta^*=\partial^t$ as operators on $L^2_tC^*(\Sigma)$.
It follows that $(\partial^t)^2=0$ (since $\delta^2=0$), and we
can consider (reduced) $L^2_t$--(co)homology:
$$L^2_tH^i(\Sigma)=\ker{\delta^i}/\overline{{\rm im }\,{\delta^{i-1}}},\quad
L^2_tH_i(\Sigma)=\ker{\partial^t_i}/\overline{{\rm im }\,{{\partial^t_{i+1}}}}
$$
Since $\delta^*=\partial^t$, $(\partial^t)^*=\delta$ we have
$L^2_tC^i(\Sigma)=\ker{\partial^t_i}\oplus\overline{{\rm im }\,{\delta^{i-1}}}
=\ker{\delta^i}\oplus\overline{{\rm im }\,{{\partial^t_{i+1}}}}$
(orthogonal direct sums).
It follows that $$L^2_tH^i(\Sigma)\simeq L^2_t{\cal H}^i(\Sigma)
\simeq L^2_tH_i(\Sigma),$$
where $ L^2_t{\cal H}^i(\Sigma)$ is the space $\ker{\delta^i}\cap
\ker{\partial^t_i}$ of harmonic $i$--cochains.
\begin{remark}
Suppose that $X(q)$ is a building associated to a $BN$--pair,
with Weyl group $W$. Then the $B$--invariant part of the 
$L^2$ cochain complex of $X(q)$ is isomorphic to 
$L^2_qC^*(\Sigma)$.
\end{remark}

\section{Hecke algebra}\label{sec:2}

We deform the usual scalar product on ${\bf C}[W]$ into $\langle\, 
,\rangle_t$:
\begin{equation}\langle \sum_{w\in W} a_w\delta_w,\sum_{w\in W} b_w\delta_w\rangle_t=
\sum_{w\in W} a_w\overline{b_w}t^{d(w)}.
\label{eq:2.1}\end{equation}
We also correspondingly deform the multiplication into the
following Hecke $t$--mul\-ti\-pli\-ca\-tion: for $w\in W$, $s\in S$
 we put
\begin{equation}\delta_w\delta_s=\begin{cases}
\delta_{ws}                &\hbox{\qquad\rm if $d(ws)>d(w)$;}\\
t\delta_{ws}+(t-1)\delta_w &\hbox{\qquad\rm if $d(ws)<d(w)$.}
\end{cases}
\label{eq:2.2}\end{equation} This extends to a ${\bf C}$--bilinear associative
multiplication on ${\bf C}[W]$ (see Bourbaki \cite{Bourb}). Using
\eqref{eq:2.2} and induction on $d(v)$ one easily shows
\begin{equation}
\delta_w\delta_v=\delta_{wv}\hbox{\qquad\rm if $d(wv)=d(w)+d(v)$.}
\label{eq:2.3}
\end{equation}
We keep the involution  on ${\bf C}[W]$ independent of $t$:
\begin{equation}
\biggl(\sum_{w\in W} a_w\delta_w\biggr)^*=
\sum_{w\in W}\overline{a_{w^{-1}}}\delta_w.\label{eq:2.4}
\end{equation}
\begin{propo}\label{prop:2.1}
The above scalar product, multiplication and involution define
a Hilbert algebra structure on ${\bf C}[W]$ (in the sense of Dixmier
\cite[A.54]{Dix});
we use the notation ${\bf C}_t[W]$ to indicate this structure.
\end{propo}

\begin{proof}
We begin with involutivity: $(xy)^*=y^*x^*$.
One checks it using \eqref{eq:2.2} and \eqref{eq:2.3} for 
$x=\delta_w$, $y=\delta_s$ considering two cases:
$d(ws)<d(w)$, $d(ws)>d(w)$. Then one checks it for
$x=\delta_w$, $y=\delta_u$ by induction on $d(u)$.
Finally, by ${\bf C}$--bilinearity of multiplication,
the result extends to general $x,y$. 
From involutivity and \eqref{eq:2.2} we immediately get
\begin{equation}
\delta_s\delta_w=\begin{cases}
\delta_{sw}                &\hbox{\qquad\rm if $d(sw)>d(w)$;}\\
t\delta_{sw}+(t-1)\delta_w &\hbox{\qquad\rm if $d(sw)<d(w)$.}
\end{cases}
\label{eq:2.5}
\end{equation}
We now recall and prove the conditions (i)--(iv) of \cite{Dix} 
defining a Hilbert algebra.

(i)\qua $\langle x,y\rangle_t=\langle y^*,x^*\rangle_t$.

This is a straightforward calculation (using $d(w)=d(w^{-1})$).

(ii)\qua $\langle xy,z\rangle_t=\langle y,x^*z\rangle_t$.

Due to linearity it is enough to check (ii) in the case $y=\delta_w$,
$z=\delta_u$, $x=\delta_v$. First one treats the case
$v=s\in S$, directly using \eqref{eq:2.5}; this requires four sub-cases,
depending on comparison of $d(sw)$ with $d(w)$ and $d(su)$ with $d(u)$.
Then one performs an easy induction on $d(v)$.

(iii)\qua For every $x\in {\bf C}_t[W]$ the map $ {\bf C}_t[W]\ni y\mapsto xy\in
{\bf C}_t[W]$ is continuous.

One checks first that $y\mapsto \delta_sy$
is continuous, directly using \eqref{eq:2.5}. 
Continuity of $y\mapsto xy$ for arbitrary $x\in {\bf C}_t[W]$ follows, 
because compositions and linear combinations of 
continuous maps are continuous. 

(iv)\qua The set $\{xy\mid x,y\in{\bf C}_t[W]\}$ is dense in 
${\bf C}_t[W]$.

This is immediate, since we have a unit element $\delta_1$ 
in ${\bf C}_t[W]$.
\end{proof}

\begin{coro}\label{cor:2.2}
The coefficient of $\delta_1$ in $ab$ is equal to $\langle a,b^*\rangle_t$.
\end{coro}
\begin{proof}
That coefficient is equal to $\langle ab,\delta_1\rangle_t$,
which by (ii) and (i) is $\langle b,a^*\rangle_t=\langle a,b^*\rangle_t$.
\end{proof}

As in \cite[A.54]{Dix}, we get two von Neumann algebras $U_t$, $V_t$: 
they are weak closures
of ${\bf C}_t[W]$ acting on its completion $L^2_t$  by left 
(respectively right) multiplication.

As in \cite[A.57]{Dix}, we put ${\bf C}_t[W]'$ to be the algebra of all bounded
elements of $L^2_t$; bounded means that left (or, equivalently, right)
multiplication by the element is bounded on ${\bf C}_t[W]$ (so, extends
to a bounded operator on $L^2_t$ and defines an element of $U_t$ or $V_t$).  
 
As in \cite[A.60]{Dix}, we have natural traces tr on $U_t$, $V_t$:
if $B\in U_t$ (or $B\in V_t$) is self-adjoint and positive, we ask whether
$B^{1\over2}=a\cdot$ 
(resp.\ $B^{1\over2}=\cdot a$) for an $a\in{\bf C}_t[W]'$.
If it is so, we put ${\rm tr}\,B=\|a\|^2_t$; otherwise we put
 ${\rm tr}\,B=+\infty$. The $a=\sum_{w\in W} a_w\delta_w$ 
we are asking for is self-adjoint: $a_w=\overline{a_{w^{-1}}}$,
so that by \fullref{cor:2.2} $\|a\|^2_t$
is equal to the coefficient of $\delta_1$ in $a^2$.
Thus $B$ is the multiplication by the bounded self-adjoint element
$b=a^2$, and ${\rm tr}\,B$ is equal to the coefficient 
 of $\delta_1$ in $b$.

Suppose now that we are given a closed subspace $Z$ of $\oplus_{i=1}^lL^2_t$,
such that the orthogonal projection $P_Z$ onto $Z$ is an element of
$M_{l\times l}\otimes V_t$. To calculate the trace of this projection
 we first need
to identify $P_Z$ as a matrix. So, we take the standard basis 
$\{e_i\}$ of $\oplus_{i=1}^lL^2_t$ ($e_i$ has $\delta_1$ as the 
$i$th coordinate,
and other coordinates 0), and apply $P_Z$ to it. We expand the results
in the basis $\{e_i\}$: let $a^j_i\in L^2_t$ be the
$j$th coordinate of $P_Z(e_i)$. Then we take the coefficient
of $\delta_1$ in $a^i_i$ and sum over $i$. The number we get is the trace
of $P_Z$.   

\section{$L^2_t$--Betti numbers}\label{sec:3}

It will be convenient to identify $L^2_t$ with $L^2(W, \nu_t)$,
where $\nu_t(w)=t^{d(w)}$. For any Coxeter group $\Gamma$
(we have $W$ as well as its subgroups $W_T$ in mind) 
 the generating function of $\Gamma$ is defined by
 $\Gamma(x)=\sum_{\gamma\in\Gamma}x^{d(\gamma)}$. 
For a finite $\Gamma$ it is a polynomial,
in general it is a rational function. We denote by $\rho_\Gamma$
the radius of convergence of the series defining $\Gamma(x)$.

As in the case of buildings (\fullref{sec:0}), we 
have 
$L^2_tC^i(\Sigma)=\bigoplus_{\sigma^i\subset D}L^2(W\sigma^i,\mu_t)$.
Now $L^2(W\sigma^i,\mu_t)$ can be identified with 
$L^2(W,\nu_t)^{W_{T(\sigma)}}$ (where $T(\sigma)$ is the largest subset
of $S$ such that $\sigma\subseteq D_{T(\sigma)}$) 
via the map $\phi$ given by $\phi(f)(w)
={1\over\sqrt{W_{T(\sigma)}(t)}}f(w\sigma)$ (we distorted the natural
map by the factor 
 ${1\over\sqrt{W_{T(\sigma)}(t)}}$ in order to make it
isometric). Finally, $L^2(W,\nu_t)^{W_{T(\sigma)}}$ is a subspace of
$L^2(W,\nu_t)=L^2_t$, so that we get an isometric embedding 
$$\Phi\colon L^2_tC^i(\Sigma)\hookrightarrow\bigoplus_{\sigma^i\subset D}L^2_t
=C^i(D)\otimes L^2_t.$$
Let ${\cal L}$ denote the algebra $U_t$ acting diagonally on the left 
on $\oplus_{\sigma\subset D}L^2_t=C^*(D)\otimes L^2_t$;
let ${\cal R}$ be ${\rm End }\,(C^*(D))\otimes V_t$ acting on the same 
space on the right. The von Neumann algebras ${\cal L}$ and ${\cal R}$
are commutants of each other.
\begin{lem}\label{lem:3.1}
The projection of $L^2_t$ onto 
$L^2(W\sigma,\mu_t)=L^2(W,\nu_t)^{W_{T(\sigma)}}$
is given by the \rhtm by 
\begin{equation}
p_{T(\sigma)}={1\over W_{T(\sigma)}(t)}
\sum_{w\in W_{T(\sigma)}}\delta_w.\label{eq:3.1}
\end{equation}
\end{lem}
\begin{proof} Put $T=T(\sigma)$.
The subspace onto which we project consists of those elements
of $L^2_t$ which are right $W_T$--invariant; this is equivalent to
being invariant  under \rhtm by ${1\over 1+t}(\delta_1+\delta_s)$ for all 
$s\in T$ (to check this one splits $W$ into pairs $\{w,ws\}$,
and calculates for each pair separately using \eqref{eq:2.2}).
As a result, this subspace is ${\cal L}$--invariant,
so that the projection $P_T$ onto it is an element of ${\cal R}$.
It follows that $P_T$ is given by \rhtm by $P_T(\delta_1)$.
The latter is clearly of the form $C\sum_{w\in W_T}\delta_w$,
where $C$ is a constant such that 
$$\langle \delta_1-C\sum_{w\in W_T}\delta_w,C\sum_{w\in W_T}\delta_w\rangle_t
=0.$$
This gives $C=\| \sum_{w\in W_T}\delta_w\|_t^{-2}=\bigl(\sum_{w\in W_T}
t^{d(w)}\bigr)^{-1}={1\over W_T(t)}$.\end{proof}
\begin{lem}\label{lem:3.2}
${\cal L}$ preserves the subspace $L^2_tC^i(\Sigma)\subset 
C^*(D)\otimes L^2_t$ and commutes with $\delta$ and $\partial^t$.
\end{lem}
\begin{proof} The first claim follows from \fullref{lem:3.1} (and actually 
was a step in the proof of that lemma).  
To prove the second part notice that $\delta$ is an element
of ${\cal R}$: the matrix with $V_t$--coefficients describing $\delta$  
has non-zero $\sigma\tau$--entry if and only if $\sigma$ is a codimension
1 face of $\tau$; the entry is then $\sqrt{W_{T(\sigma)}(t)\over
W_{T(\tau)}(t)}\delta_1$. It follows that $\delta$ commutes with ${\cal L}$.
So therefore does its adjoint $\partial^t$. \end{proof}
\begin{coro}\label{cor:3.3}
$L^2_tC^i(\Sigma)$, $L^2_t{\cal H}^i(\Sigma)$, $\ker{\delta^i}$,
$\ker{\partial^t_i}$, $\overline{{\rm im }\,\delta^i}$,
$\overline{{\rm im }\,\partial^i_t}$ are ${\cal L}$--in\-va\-riant;
therefore, orthogonal projections onto these spaces
belong to ${\cal R}$.
\end{coro}
We use tr to denote the tensor product of the usual matrix 
trace on ${\rm End }\,(C^*(D))$ and the von Neumann trace on $V_t$
as described in \fullref{sec:2}. We put
\begin{gather}
b^i_t=L^2_tb^i(\Sigma)={\rm tr}\left(\hbox{\rm projection onto }
L^2_t{\cal H}^i(\Sigma)\right)\label{eq:3.2}\\
c^i_t=L^2_tc^i(\Sigma)={\rm tr}\left(\hbox{\rm projection onto }
L^2_tC^i(\Sigma)\right)\label{eq:3.3}\\
\chi_t=\sum_i(-1)^ib^i_t=\sum_i(-1)^ic^i_t\,.\label{eq:3.4}\end{gather}
The sums in \eqref{eq:3.4} give the same value by the standard 
algebraic topology argument. It follows from \fullref{lem:3.1}
that $c^i_t=\sum_{\sigma^i\subset D}{1\over W_{T(\sigma)}(t)}$.
Grouping together simplices $\sigma$ with the same $T(\sigma)$ and using
formula (5) from Charney--Davis \cite{ChD} we obtain the following result 
(see Serre \cite{Se}).
\begin{coro}\label{cor:3.4}
$$\chi_t={1\over W(t)}$$
\end{coro}
\begin{thm}\label{thm:3.5}
Suppose that $X(q)$ is a building associated to a $BN$--pair,
with Weyl group $W$. Then $L^2b^i(X(q))=b^i_q$.
\end{thm}
\begin{proof} For $t=q$, $L^2_tC^i(\Sigma)$ coincides with the space of 
$B$--invariant elements of $L^2C^i(X(q))$.
By the concluding remarks of \fullref{sec:0}, the matrix of the projection
onto $L^2{\cal H}^i(X(q))$ has $B$--invariant entries---so that
it coincides with the one we use to define $b^i_t$. Hence the conclusion.
\end{proof}
Suppose now that the pair $(D,\partial D=D\cap\partial\Delta)$ 
is a generalised homology $n$--disc (ie, it is a homology manifold with 
boundary, with relative homology groups the same as those of an $n$--disc
modulo its boundary). Then each $D_T=D\cap\Delta_T$ is also a 
homology $(n-|T|)$--disc (for $T\in{\cal F}$). We can now use 
$wD_T$, $w\in W$, $T\in{\cal F}$, as a homology cellular structure
on $\Sigma$ (denoted $\Sigma_{ghd}$). The cell $D_T$ has the form 
of an $o_T$--centred cone; we put $d(wD_T)=d(wo_T)$, and define 
$\mu_t$, (co)chain complexes, the embedding $\Phi$, the $U_t$--module
structure and the numbers $b^i_t(\Sigma_{ghd})$ in essentially the same way as 
for the original triangulation of $\Sigma$.


\section{Dual cells}\label{sec:4}

So far we used the triangulation of $\Sigma$ which originated from
the barycentric subdivision of a simplex. We will use notation 
$\Sigma_{st}$ to remind that we have this standard triangulation in mind.
In this section we will describe another cell structure on $\Sigma$.
It will
make our discussion of Poincar\' e duality in \fullref{sec:6} 
look pretty standard.

To each $T\in{\cal F}$ we associate a
face $\Delta_T$ of $\Delta$, whose barycentre $o_T$ is a vertex
of the Davis chamber $D$. We define $\langle T\rangle$ as the union
of all simplices $\sigma\subset\Sigma$ such that $\sigma\cap D_T=o_T$
(recall that $D_T=D\cap\Delta_T$). As a simplicial complex,
$\langle T\rangle$ is an $o_T$--centred cone over $\Sigma_T$; since $T$ is such that $W_T$ 
is finite, $\Sigma_T$ is a sphere and $\langle T\rangle$ is a disc
of dimension $|T|$. 
The boundary of $\langle T\rangle$ is cellulated by
$w\langle U\rangle$, for all possible $T\subset U\subseteq S$, $w\in W_T$.   
The complex $\Sigma$ cellulated by $w\langle T\rangle$, over all 
$w\in W$, $T\in\cal{F}$, is a cellular complex 
that we denote $\Sigma_d$. 
The cells of $\Sigma_d$ will be called \it dual cells\rm .
The name \it Coxeter blocks \rm is also used (Davis \cite{D1}).


We now put $d(w\langle T\rangle)=d(wo_T)$, and define the measures $\mu_t$ 
on the set $\Sigma_d^{(i)}$
of $i$--dimensional cells of $\Sigma_d$ by $\mu_t(\langle a\rangle)=
t^{d(\langle a\rangle)}$.
Then
$$L^2_tC^i(\Sigma_d)=L^2_tC_i(\Sigma_d)\simeq L^2(\Sigma_d^{(i)},\mu_t),$$
We now define $\delta^i\colon L^2_tC^i(\Sigma_d)\to L^2_tC^{i+1}(\Sigma_d)$ by
$$\delta^i(f)(\langle \tau\rangle^{i+1})=
\sum_{\langle \sigma\rangle^i\subset\langle \tau\rangle^{i+1}}
[\langle \tau\rangle:\langle \sigma\rangle]f(\langle \sigma\rangle)$$
and  $\partial_i^t\colon L^2_tC_i(\Sigma_d)\to L^2_tC_{i-1}(\Sigma_d)$ by
$$\partial_i^t(f)(\langle \eta\rangle^{i-1})=\sum_{\langle \sigma\rangle^i\supset\langle \eta\rangle^{i-1}}
[\langle \eta\rangle:\langle \sigma\rangle]t^{d(\langle \sigma\rangle)-d(\langle \eta\rangle)}f(\langle \sigma\rangle).$$
The discussion from \fullref{sec:1} can be continued,
and supplies us with $L^2_t{\cal H}^i(\Sigma_d)$.
Now we wish to bring in the Hecke algebra.
We pick (arbitrarily) orientations of the cells $\langle T\rangle$
($T\in{\cal F}$), and extend these to orientations of all
cells in $\Sigma_d$ as follows: $w\langle T\rangle$ is the oriented cell
which is the image of the oriented cell $\ce{T}$ by $w$,
with orientation changed by a factor of $(-1)^{d(w)}$.    
Using these orientations, we identify $L^2_tC^*(\Sigma_d)$
with $\oplus_{T\in{\cal F}}L^2(W\ce{T},\mu_t)$.
For every $T\in{\cal F}$ we define a map $\psi_T\colon L^2(W\langle T\rangle,\mu_t)
\to L^2_t$ by the formula
\begin{equation}
\psi_T(f)=\sum_{w\in W^T} f(w\ce{T}) 
(-1)^{d(w)}\sqrt{W_T(t^{-1})}\,\delta_wh_T,
\label{eq:4.1}
\end{equation}
where $W^T=\{w\in W\mid {\forall u\!\in\! W_T}, d(wu)\geq d(w)\}$
(the set of  $T$--reduced  elements),
and 
\begin{equation}
h_T={1\over W_T(t^{-1})}\sum_{u\in W_T}(-t)^{-d(u)}\delta_u.\label{eq:4.2}
\end{equation}
Putting together these maps we get a map $\Psi\colon L^2_tC^*(\Sigma_d)
\to\oplus_{T\in{\cal F}}L^2_t$.

\begin{lem}\label{lem:4.1}
{\rm(1)}\qua For all $s\in T$ we have $\delta_sh_T=-h_T$.

{\rm(2)}\qua For all $u\in W_T$ we have $\delta_uh_T=(-1)^{d(u)}h_T$.

{\rm(3)}\qua For all $U\subseteq T$ we have $h_Uh_T=h_T$. 
\end{lem}
\begin{proof} (1)\qua Let $w\in W$ be such that $d(sw)>d(w)$. Then $\delta_s\delta_w=
\delta_{sw}$ (by \eqref{eq:2.3}).
We then have
{\setlength\arraycolsep{2pt}
\begin{eqnarray}
\delta_s(\delta_w-{1\over t}\delta_{sw})
&=&\delta_{sw}-{1\over t}(\delta_s\delta_s)\delta_w
=\delta_{sw}-{1\over t}
(t\delta_1+(t-1)\delta_s)\delta_w\nonumber\\
&=&(1-{t-1\over t})\delta_{sw}-\delta_w=-(\delta_w-{1\over t}\delta_{sw})\nonumber
\end{eqnarray}
}%
Since $h_T$ is a linear combination of expressions of the form
$\delta_w-{1\over t}\delta_{sw}$, (1) follows.  

(2)\qua Follows from (1) by induction on $d(u)$.

(3)\qua $\displaystyle{h_Uh_T
={1\over W_U(t^{-1})}\sum_{u\in W_U}(-t)^{-d(u)}\delta_uh_T}$\newline
$\displaystyle{\phantom{{\rm (3)\qua } h_Uh_T}={1\over W_U(t^{-1})}\sum_{u\in W_U}(-t)^{-d(u)}(-1)^{d(u)}h_T}$\newline
$\displaystyle{\phantom{{\rm (3)\qua } h_Uh_T}={1\over W_U(t^{-1})}\left(\sum_{u\in W_U}t^{-d(u)}\right)h_T=
h_T}$
\end{proof}
\begin{lem}\label{lem:4.2}
{\rm(1)}\qua For every $T\in {\cal F}$ the map $\psi_T$ is an isometric
embedding.\par
{\rm(2)}\qua The orthogonal projection of $L^2_t$ onto the image of $\psi_T$ is
given by \rhtm by $h_T$.
\end{lem}
\begin{proof} (1)\qua The squared norm of a summand from the right hand side
of \eqref{eq:4.1} is 
$$\|f(w\ce{T})(-1)^{d(w)}\sqrt{W_T(t^{-1})}\,\delta_wh_T\|^2_t=
|f(w\ce{T})|^2W_T(t^{-1})\|\delta_wh_T\|^2_t.
$$ 
Since $w$ is $T$--reduced, we have $\delta_w\delta_u=\delta_{wu}$
for all $u\in W_T$. Therefore
{\setlength\arraycolsep{2pt}
\begin{eqnarray}
\|\delta_wh_T\|^2_t&=&
\left\|{1\over W_T(t^{-1})}\sum_{u\in W_T}(-t)^{-d(u)}\delta_{wu}\right\|^2_t
\hskip-2pt=
\left|{1\over W_T(t^{-1})}\right|^2\!\sum_{u\in W_T}\!|-t|^{-2d(u)}t^{d(wu)}\nonumber\\
&=&t^{d(w)}{1\over W_T(t^{-1})^2}\!\sum_{u\in W_T}\!t^{-d(u)}=
t^{d(w)}{1\over W_T(t^{-1})}.\nonumber
\end{eqnarray}
}

(2)\qua Due to $h_Th_T=h_T$ and $h_T^*=h_T$, \rhtm by $h_T$ is an orthogonal projection.
Let $w\in W$; write $w=vu$ where $u\in W_T$ and $v$ is $T$--reduced.
Then $\delta_wh_T=\delta_v\delta_uh_T=(-1)^{d(u)}\delta_vh_T$.
This shows that image 
of the space of finitely
supported functions (on $W\ce{T}$) under $\psi_T$ is equal to the image 
of the space of finitely supported functions (on $W$) under \rhtm by $h_T$. 
Since $\psi_T$ is isometric, the $L^2_t$--completions of these images also coincide.    
\end{proof}
Denote by $\cal L$ the algebra $U_t$ acting diagonally
on the left on $\oplus_{T\in{\cal F}}L^2_t$, and by $\cal R$ its commutant
$M_{|\cal F|}({\bf C})\otimes V_t$ (acting on the right).
It follows from \fullref{lem:4.2} that the image of $\Psi$ is $\cal L$--invariant.
In other words, we have a $U_t$--module structure
on $L^2_tC^*(\Sigma_d)$, defined by the condition that the isometric 
embedding  $\Psi\colon L^2_tC^*(\Sigma_d)\to\oplus L^2_t$ is a morphism of 
$U_t$--modules. 
Thus, we think of $L^2_tC^*(\Sigma_d)$ as of a submodule
of $\oplus_{T\in{\cal F}}L^2_t$. 
\begin{lem}\label{lem:4.3}
The map $\delta\colon L^2_tC^*(\Sigma_d)\to L^2_tC^*(\Sigma_d)$ is (a restriction of)
an element of $\cal R$. For $U\subset T\in{\cal F}$ satisfying 
$|T|=|U|+1$, the $UT$--entry of this element is
$$[\ce{T}:\ce{U}]\sqrt{W_T(t^{-1})\over W_U(t^{-1})}\,h_T$$
\end{lem}
\begin{proof} Consider a pair of cells $w\ce{U}$, $w\ce{T}$.
We have $[w\ce{T}:w\ce{U}]=[\ce{T}:\ce{U}]$. 
We can assume that $w$ is $U$--reduced, and write it as 
$vu$, where $v$ is $T$--reduced
and $u\in W_T$. 
Let $f\in L^2_tC^{\dim\ce{U}}(\Sigma_d)$. The summand in 
$\psi_U(f)$ corresponding to the cell $w\ce{U}$ is
$$f(w\ce{U})(-1)^{d(w)}\sqrt{W_U(t^{-1})}\,\delta_wh_U.$$
The summand in $\psi_T(\delta f)$ corresponding to
the contribution of $f(w\ce{U})$ to\break $(\delta f)(w\ce{T})$ is
$$[\ce{T}:\ce{U}]f(w\ce{U})(-1)^{d(v)}\sqrt{W_T(t^{-1})}\,\delta_vh_T.$$
Now $\delta_wh_Uh_T=\delta_wh_T=\delta_v\delta_uh_T=(-1)^{d(u)}\delta_vh_T$,
and the lemma follows.
\end{proof}
\begin{coro}\label{cor:4.4}
The 
subspaces $L^2_tC^i(\Sigma_d)$, $L^2_t{\cal H}^i(\Sigma_d)$,
$\ker{\delta^i}$,
$\ker{\partial^t_i}$, 
$\overline{{\rm im }\,\delta^i}$ and
$\overline{{\rm im }\,\partial^i_t}$ of $\oplus_{T\in{\cal F}}L^2_t$
are ${\cal L}$--in\-va\-riant; 
therefore, orthogonal projections onto these spaces
are elements of ${\cal R}$.
\end{coro}

\section{Invariance}\label{sec:5}

In this section we prove that $L^2_tH^*(\Sigma_d)\simeq
L^2_tH^*(\Sigma_{st})$ ($\simeq L^2_tH^*(\Sigma_{ghd})$, if the latter exists)
as $U_t$--modules. It will be convenient for us to work with homology
rather than cohomology; since both are isomorphic to the $U_t$--module of harmonic
cochains, it makes no difference.  

We start by fixing orientation conventions. Let us pick arbitrary orientations
of the dual cells $\ce{T}$ for all $T\in{\cal F}$. We extend these orientations to
all dual cells as in \fullref{sec:4} ($w\ce{T}$ is oriented by $(-1)^{d(w)}$ times the orientation 
of $\ce{T}$ pushed forward by $w$).
For $T\in {\cal F}$ of cardinality $k$, let $\ce{T}\cap D^{(k)}$ be the set of
all $k$--simplices of $\Sigma_{st}$ contained in $\ce{T}\cap D$. We orient
every element of  $\ce{T}\cap D^{(k)}$ by the restriction of the chosen orientation
of $\ce{T}$. We then extend these orientations $W\!$--equivariantly (to a part 
of $\Sigma_{st}$), and put arbitrary equivariant orientations on the rest of
$\Sigma_{st}$. Notice that if a $k$--simplex $\sigma$ is contained in $w\ce{T}$
(where $T$ has cardinality $k$), then the orientation of $\sigma$
agrees with $(-1)^{d(\sigma)}$ times that of $w\ce{T}$.
Orientations being chosen, we treat (co)chains as functions on the set
of cells/simplices.

We define a topological embedding of Hilbert spaces
$\theta\colon L^2_tC^*(\Sigma_d)\to L^2_tC^*(\Sigma_{st})$.
\begin{definition}
Let $f\in L^2_tC^k(\Sigma_d)$, $\sigma\in\Sigma_{st}^{(k)}$.\par
(1)\qua If there exists $\ce{\alpha}\in\Sigma_d^{(k)}$ such that $\sigma\subseteq\ce{\alpha}$
(there is at most one such $\ce{\alpha}$), then
$$\theta\! f(\sigma)=(-1)^{d(\sigma)}t^{d(\ce{\alpha})-d(\sigma)}f(\ce{\alpha}).$$    

(2)\qua If there is no $\ce{\alpha}$ as in (1), we put 
$\theta\! f(\sigma)=0$.
\end{definition}

\begin{lem}\label{lem:5.1}
$$\partial^t\theta=\theta\partial^t.$$  
\end{lem}
\proof We will show that for all $f\in L^2_tC^k(\Sigma_d)$, $\sigma\in\Sigma_{st}^{(k)}$ we have
$\partial^t\theta\! f(\sigma)=\theta\partial^t\! f (\sigma)$. There are two cases 
to consider.

(1)\qua Suppose that there exists $\ce{\alpha}\in\Sigma_d^{(k)}$ such that $\sigma\subseteq\ce{\alpha}$.
Then
{\setlength\arraycolsep{2pt}
\begin{eqnarray}
\theta\partial^t\! f(\sigma)&=&(-1)^{d(\sigma)}t^{d(\ce{\alpha})-d(\sigma)}\partial^t\! f(\ce{\alpha})\phantom{\left|{1\over T}\right|}\nonumber\\
&=&(-1)^{d(\sigma)}t^{d(\ce{\alpha})-d(\sigma)}\sum_{\ce{\beta}^{k+1}\supset\ce{\alpha}}[\ce{\beta}:\ce{\alpha}]
t^{d(\ce{\beta})-d(\ce{\alpha})}f(\ce{\beta})\nonumber\\
&=&(-1)^{d(\sigma)}\sum_{\ce{\beta}^{k+1}\supset\ce{\alpha}}[\ce{\beta}:\ce{\alpha}]
t^{d(\ce{\beta})-d(\sigma)}f(\ce{\beta}).\label{eq:5.1}
\end{eqnarray}
}%
On the other hand,
\begin{equation}
\partial^t\theta\! f(\sigma)=\sum_{\tau^{k+1}\supset\sigma}[\tau:\sigma]t^{d(\tau)-d(\sigma)}
\theta\! f(\tau).\label{eq:5.2}
\end{equation}
Notice that if  $\theta\! f(\tau)\ne 0$ then there exists a 
dual cell $\ce{\beta}^{k+1}\supset\tau$.
Such $\ce{\beta}$ is unique and $\ce{\tau}$ is the only $(k+1)$--simplex in 
$\ce{\beta}$ with face $\ce{\sigma}$.
Therefore \eqref{eq:5.2} equals
{\setlength\arraycolsep{2pt}
\begin{eqnarray}
\sum_{\ce{\beta}^{k+1}\supset\ce{\alpha}}&&\hbox{\hskip-.4cm}[\tau:\sigma]t^{d(\tau)-d(\sigma)}
(-1)^{d(\tau)}t^{d(\ce{\beta})-d(\tau)}f(\ce{\beta})\nonumber\\
&&\hbox{\hskip-1.15cm}=\sum_{\ce{\beta}^{k+1}\supset\ce{\alpha}}[\tau:\sigma]
(-1)^{d(\tau)}t^{d(\ce{\beta})-d(\sigma)}f(\ce{\beta}).\label{eq:5.3}
\end{eqnarray}
}%
Now \eqref{eq:5.3} and \eqref{eq:5.1} are equal because
$[\tau:\sigma]=(-1)^{d(\tau)}(-1)^{d(\sigma)}[\ce{\beta}:\ce{\alpha}]$.

(2)\qua The smallest dual cell $\ce{\alpha}$ containing $\sigma$ is of dimension $m>k$.
Then $\theta\partial^t\! f(\sigma)=0$. On the other hand,
$$\partial^t\theta\! f(\sigma)= \sum_{\tau^{k+1}\supset\sigma}[\tau:\sigma]    
t^{d(\tau)-d(\sigma)}\theta\! f(\tau).$$  
Let $\tau^{k+1}\supset\sigma$, and let $\ce{\beta}\supset\ce{\alpha}$ be the smallest
dual cell containing $\tau$. If $\theta\! f(\tau)\ne 0$, then $\dim{\ce{\beta}}=k+1$,
which forces $\ce{\beta}=\ce{\alpha}$ and $m=\dim{\ce{\alpha}}=k+1$.
Thus, we are reduced to the case $m=k+1$. In this case, there are 
exactly two simplices $\sigma_{\pm}\in\Sigma_{st}^{(k+1)}$, 
$\sigma_{\pm}\subset\ce{\alpha}$, $\sigma_{\pm}\supset\sigma$.
Since $\sigma_{\pm}$ is oriented by $(-1)^{d(\sigma_{\pm})}$ times
the orientation of $\ce{\alpha}$, we have
\begin{equation}
(-1)^{d(\sigma_+)}[\sigma_+:\sigma]=
-(-1)^{d(\sigma_-)}[\sigma_-:\sigma].\label{eq:5.4}
\end{equation}
Therefore
{\setlength\arraycolsep{2pt}
\begin{eqnarray}
\hspace{1.5cm}\partial^t\theta\! f(\sigma)&=&[\sigma_+:\sigma]t^{d(\sigma_+)-d(\sigma)}\theta\! f(\sigma_+)
+[\sigma_-:\sigma]t^{d(\sigma_-)-d(\sigma)}\theta\! f(\sigma_-)\nonumber\\
&=&[\sigma_+:\sigma]t^{d(\sigma_+)-d(\sigma)}(-1)^{d(\sigma_+)}t^{d(\ce{\alpha})-d(\sigma_+)}
f(\ce{\alpha})\nonumber\\
&&+[\sigma_-:\sigma]t^{d(\sigma_-)-d(\sigma)}(-1)^{d(\sigma_-)}t^{d(\ce{\alpha})-d(\sigma_-)}
f(\ce{\alpha})\nonumber\\
&=&((-1)^{d(\sigma_+)}[\sigma_+:\sigma]+
(-1)^{d(\sigma_-)}[\sigma_-:\sigma])t^{d(\ce{\alpha})}f(\ce{\alpha})\nonumber\\ 
&=&0.\hspace{9.4cm}\sq
\end{eqnarray}
}%
\begin{lem}\label{lem:5.2}
$\theta$ is a morphism of $U_t$--modules.
\end{lem}
\begin{proof} The $U_t$--module structures on $L^2_tC^k(\Sigma_d)$ and on $L^2_tC^k(\Sigma_{st})$
 are defined via embeddings $\Psi$ and $\Phi$. We will compare $\Psi$ and $\Phi\circ\theta$.
Let $f\in L^2_tC^k(\Sigma_d)$; $\Psi(f)$ is a collection of $\psi_T(f)$, where
\begin{equation}
\psi_T(f)=\sum_{w\in W^T}f(w\ce{T})(-1)^{d(w)}\sqrt{W_T(t^{-1})}\,\delta_wh_T.\label{eq:5.5}
\end{equation}
The part of $\theta\! f$ corresponding to $\psi_T(f)$ is 
supported by the set of $W\!$--translates 
of  simplices $\sigma\in\ce{T}\cap D^{(k)}$, and is mapped by $\Phi$ into
$\oplus_{\sigma\in\ce{T}\cap D^{(k)}}L^2_t$. The component indexed by $\sigma$ is 
$\sum_{w\in W}\theta\!f(w\sigma)\delta_w$ (notice that the stabiliser of $\sigma$ is
trivial), ie,
\begin{equation}
\sum_{w\in W}(-1)^{d(w\ce{T})}t^{d(w\ce{T})-d(w\sigma)}f(w\ce{T})\delta_w.
\label{eq:5.6}
\end{equation} 
Comparing \eqref{eq:5.5} and \eqref{eq:5.6} with the help of \eqref{eq:4.2}, we get that $\psi_T(f)$ agrees
with (every component of) the corresponding part of $\Phi(\theta\! f)$,
up to a multiplicative factor of $\sqrt{W_T(t^{-1})}$. This implies the lemma.\end{proof}

\begin{thm}\label{thm:5.3}
The map $\theta$ induces an isomorphism of $U_t$--modules $L^2_tH_*(\Sigma_d)\simeq
L^2_tH_*(\Sigma_{st})$.  
\end{thm}
\begin{proof} Lemmas \ref{lem:5.1} and \ref{lem:5.2} imply that $\theta$ induces a morphism of
$U_t$--modules on homology. We have to check that it is an isomorphism of vector spaces.

Let $K_*$ be the image of $\theta$. It is a subcomplex of $(L^2_tC_*(\Sigma_{st}),\partial^t)$.
A $k$--chain $c\in L^2_tC_*(\Sigma_{st})$ is in $K_*$ if and only if the following
two conditions hold:

(1)\qua $c$ is supported by the union of $k$--dimensional 
dual cells: $\bigcup\Sigma_{d}^{(k)}$;

(2)\qua if $\sigma^k,\tau^k\subseteq\ce{\alpha}^k$, then $c(\sigma)=(-t)^{d(\tau)-d(\sigma)}c(\tau)$.

We need to show that the inclusion 
$K_*\hookrightarrow L^2_tC_*(\Sigma_{st})$ 
induces an isomorphism on (reduced) homology.

Let $m_t\colon L^2_tC_*(\Sigma_{st})\to L^2_{t^{-1}}C_*(\Sigma_{st})$ be the isomorphism
(of Hilbert spaces) $m_tf(\sigma)$ $=t^{d(\sigma)}f(\sigma)$. Instead of 
working directly with $K_*$, $L^2_tC_*(\Sigma_{st})$ and $\partial^t$, we will work with
$L_*=m_t(K_*)$, $E_*=L^2_{t^{-1}}C_*(\Sigma_{st})=m_t(L^2_tC_*(\Sigma_{st}))$ and
$\partial=m_t\partial^t m_t^{-1}$. The advantage is that
{\setlength\arraycolsep{2pt}
\begin{eqnarray}
\partial g(\sigma)&=&
m_t\partial^t m_t^{-1}g(\sigma)=
t^{d(\sigma)}\partial^t m_t^{-1}g(\sigma)\phantom{\left|{1\over T}\right|}
\nonumber\\
&=&
t^{d(\sigma)}\sum_{\tau^{k+1}\supset\sigma} [\tau:\sigma]t^{d(\tau)-d(\sigma)}m_t^{-1}g(\tau)\hbox{\hskip1.7cm}\label{eq:5.7}\\
&=&
\sum_{\tau^{k+1}\supset\sigma} [\tau:\sigma]t^{d(\tau)}t^{-d(\tau)}g(\tau)=
\sum_{\tau^{k+1}\supset\sigma} [\tau:\sigma]g(\tau).\hbox{\hskip1.7cm}\nonumber
\end{eqnarray}
}%
To check whether $c\in E_*$ is in $L_*$ we use (1) and the following version of (2):

$(2')$\qua  if $\sigma^k,\tau^k\subseteq\ce{\alpha}^k$, then $c(\sigma)=(-1)^{d(\tau)-d(\sigma)}c(\tau)$.
\begin{lem}\label{lem:5.4}
Let $c\in E_k$. If $\partial c\in L_*$, then there exists
a $d\in E_{k+1}$ such that $c-\partial d\in L_*$. Moreover, there is a constant $C$
depending only on $W$ and $t$ such that $d$ can be chosen so that $\|d\|\le C\|c\|$.
\end{lem}
\begin{proof} Each dual cell $\ce{\alpha}$ is a disc; we denote by ${\rm int}\ce{\alpha}$ its interior,
and by ${\rm bd}\ce{\alpha}$ its boundary.
We construct, by descending induction on $m$ ($m\geq k$), cochains $d_m\in E_{k+1}$
such that $c-\partial d_m$ is supported by the union of dual cells of dimensions at most $m$. 
For $m\geq\dim{\Sigma}$ we put $d_m=0$. Suppose that $d_m$ is already constructed, where $m>k$.
For every dual $m$--cell $\ce{\alpha}$, let $c_\alpha$ be the restriction of 
$c-\partial d_m$ to $\ce{\alpha}$ (ie, if   $c-\partial d_m=\sum a_\sigma\sigma$, then $c_\alpha=
\sum_{\sigma\subseteq\ce{\alpha}}a_\sigma\sigma$). 
Let $\sigma^k\cap{\rm int}\ce{\alpha}\ne\emptyset$. Then $\sigma$ appears in $\partial c_\alpha$
and in $\partial c=\partial (c-\partial d)$ with the same coefficient,
due to the inductive assumption. But, since $\partial c\in L_*$, this coefficient is 0.
As a result, $c_\alpha\in Z_k(\ce{\alpha},{\rm bd}\ce{\alpha})$. 
Since $H_k(\ce{\alpha},{\rm bd}\ce{\alpha})=0$ (recall that $m=\dim{\ce{\alpha}}>k$),
we can find $d_\alpha\in C_{k+1}(\ce{\alpha})$ such that $c_\alpha-\partial d_\alpha\in C_k({\rm bd}
\ce{\alpha})$. Moreover, we can choose $d_\alpha$ so that $\|d_\alpha\|\leq C_1\|c_\alpha\|$, 
for some constant $C_1$ depending only on $W$ and $t$. Due to uniform local finiteness
of $\Sigma$, we deduce $\|\sum_{\ce{\alpha}} d_\alpha\|\leq C_2\|c\|$ for some constant $C_2$.
We put $d_{m-1}=d_m+\sum_{\ce{\alpha}\in\Sigma_d^{(m)}}d_\alpha$, and $d=d_k$.        
 
The estimate $\|d\|\leq C\|c\|$ clearly follows from the construction.
The chain $c-\partial d=\sum b_\sigma\sigma$ is supported by the union
of dual cells of dimensions at most $k$. Let us check that it satisfies the condition $(2')$.
Suppose that $\sigma^{k-1}\cap{\rm int}\ce{\alpha}^k\ne\emptyset$. There are exactly two $k$--simplices
$\sigma_{\pm}\subset\ce{\alpha}$ such that $\sigma\subset\sigma_{\pm}$. The coefficient
of $\sigma$ in $\partial (c-\partial d)=\partial c$ is 0 (because $\partial c\in L_*$),
and, on the other hand, is equal to $[\sigma_+:\sigma]b_{\sigma_+}+
[\sigma_-:\sigma]b_{\sigma_-}$. Using \eqref{eq:5.4} we get $b_{\sigma_+}=
(-1)^{d(\sigma_+)-d(\sigma_-)}b_{\sigma_-}$.
This holds for all $\sigma^{k-1}$ satisfying $\sigma^{k-1}\cap{\rm int}\ce{\alpha}^k\ne\emptyset$, 
which implies that
$c-\partial d$ satisfies $(2')$. Hence $c-\partial d\in L_*$. The lemma is proved.
\end{proof} 
We are ready to check that the inclusion $\iota\colon L_*\hookrightarrow E_*$ induces
an isomorphism $\iota_*$ on  (reduced) homology. To show that $\iota_*$ is surjective,
suppose that $c\in E_*$ is closed: $\partial c=0$. Then $\partial c\in L_*$,
and, by \fullref{lem:5.4}, there exists $d\in E_*$ such that $c-\partial d\in L_*$.
We get $[c]=\iota_*[c-\partial d]$.

To show that $\iota_*$ is 1--1, suppose that $l\in L_*$, $\partial l=0$ and
$\iota_*[l]=0$, ie,  
$l=\lim{\partial e_n}$ for some sequence of $e_n\in E_*$. Applying \fullref{lem:5.4}
to $c=l-\partial e_n$, we get that there exist $f_n\in E_*$, $f_n\to 0$ such that
$l-\partial e_n-\partial f_n\in L_*$. But, since $l\in L_*$,
we deduce that $\partial (e_n+f_n)\in L_*$. Now we apply \fullref{lem:5.4} to
$c=e_n+f_n$ to get $g_n\in E_*$ such that $h_n=e_n+f_n-\partial g_n\in L_*$.
We have $$\partial h_n=\partial e_n+\partial f_n-\partial \partial g_n.$$ 
The last term is 0, the middle term converges to 0 since
$\partial $ is bounded and $f_n\to0$, so that, finally,
$$\lim{\partial h_n}=\lim{\partial e_n}=l.$$ 
This means that $[l]=0$ in $H_*(L_*)$.

We have shown that $(L_*,\partial )\hookrightarrow(E_*,\partial )$ induces an isomorphism 
on homology. Therefore so does the inclusion $(K_*,\partial^t)\hookrightarrow(L^2_tC_*(\Sigma_{st}),\partial^t)$.
The theorem follows.
\end{proof}
Let us now assume that $D$ is a generalised homology disc. Then, along the same lines as above,
one shows $L^2_tH^*(\Sigma_{st})\simeq L^2_tH^*(\Sigma_{ghd})$ (as $U_t$--modules).
More precisely, one defines $\theta\colon L^2_tH^*(\Sigma_{ghd})\to L^2_tH^*(\Sigma_{st})$
by $\theta\! f(\sigma)=f(\alpha)$ if $\sigma^k\subseteq\alpha^k\in\Sigma^{(k)}_{ghd}$,
and $\theta\! f(\sigma)=0$ if no such $\alpha^k$ exists. The proof of $\partial^t\theta=\theta\partial^t$
is similar to that of \fullref{lem:5.1}, and it is clear that $\theta$ is a $U_t$--morphism.
A chain $c\in L^2_tC_k(\Sigma_{st})$ is in the image $K_*$ of $\theta$ 
if and only if

(1)\qua $c$ is supported by $\bigcup \Sigma_{ghd}^{(k)}$;

(2)\qua if $\sigma^k,\tau^k\subseteq\alpha^k\in\Sigma_{ghd}^{(k)}$, then 
$c(\sigma)=c(\tau)$.

These conditions do not change under $m_t$, and the rest of the proof of 
\fullref{thm:5.3}
can be repeated with dual cells replaced by cells of $\Sigma_{ghd}$
(the only other change will be
$[\sigma_+:\sigma]=-[\sigma_-:\sigma]$ instead of the more complicated \eqref{eq:5.4}).
We get
\begin{thm}\label{thm:5.5}
Let $(D,\partial D)$ be a generalised homology disc. 
Then we have the following isomorphisms of (graded) $U_t$--modules: 
$L^2_tH^*(\Sigma_{ghd})\simeq L^2_tH^*(\Sigma_{st})\simeq 
L^2_tH^*(\Sigma_d)$.    
\end{thm}

\section{Poincar\' e Duality}\label{sec:6}

Let us define a map $D\colon L^2_t\to L^2_{t^{-1}}$
by 
\begin{equation}
D(\sum a_w\delta_w)=\sum (-t)^{d(w)}a_w\delta_w.\label{eq:6.1}
\end{equation}
Direct calculation shows that $D$ is an isometric isomorphism
of Hilbert spaces.
Notice that $D$ maps ${\bf C}_t[W]$ onto ${\bf C}_{t^{-1}}[W]$.   
It is easy to check that $D$ preserves the relations 
defining Hecke multiplication:
if $d(ws)>d(w)$, then $$D(\delta_w\delta_s)=D(\delta_{ws})=
(-t)^{d(ws)}\delta_{ws}=(-t)^{d(w)}\delta_w(-t\delta_s)=
D(\delta_w)D(\delta_s);$$
if $d(ws)<d(w)$, then 
{\setlength\arraycolsep{2pt}
\begin{eqnarray}
D(\delta_w\delta_s)&=&
D(t\delta_{ws}+(t-1)\delta_w)=
t(-t)^{d(ws)}\delta_{ws}+(t-1)(-t)^{d(w)}\delta_w\nonumber\\
&=&
(-t)^{d(w)+1}t^{-1}\delta_{ws}+(-t)^{d(w)+1}(t^{-1}-1)\delta_w=
(-t)^{d(w)}\delta_w(-t)\delta_s\nonumber\\
&=&
D(\delta_w)D(\delta_s).\nonumber
\end{eqnarray}
}%
Hence, $D$ restricts to an isometric isomorphism of Hilbert 
algebras
${\bf C}_t[W]$ and ${\bf C}_{t^{-1}}[W]$. In particular, $D$ preserves
products:
for all $x,y\in {\bf C}_t[W]$, we have $D(xy)=D(x)D(y)$.   
Passing to limits with $y$ in the norm $\|\cdot\|_t$, we deduce
that  
the map $D\colon L^2_t\to L^2_{t^{-1}}$ is a morphism of left modules
over the algebra morphism $D\colon {\bf C}_t[W]\to{\bf C}_{t^{-1}}[W]$.
Then passing to limits with $x$ in the weak operator topology,
we deduce that $D\colon L^2_t\to L^2_{t^{-1}}$ is a morphism of left modules
over the von Neumann algebra isomorphism $D\colon U_t\to U_{t^{-1}}$.
Analogous statements hold for the right module structures.
Finally, since $D$ preserves the coefficient of
$\delta_1$, it preserves dimensions of (left) submodules
of $L^2_t$.
\begin{thm}\label{thm:6.1}
Suppose that the pair $(D,\partial D)$ 
is a generalised homology $n$--disc. Then $b^i_t=b^{n-i}_{t^{-1}}$.
\end{thm}
\begin{proof} There is a bijection $D_T\leftrightarrow\langle T\rangle$, 
where $T\in{\cal F}$;
it can be unambiguously extended to $wD_T\leftrightarrow w\langle T\rangle$,
a natural bijection between $i$--cells of $\Sigma_{ghd}$ and
$(n-i)$--cells of $\Sigma_d$. When $w$ and $T$ are not specified we 
write simply $\sigma\leftrightarrow\langle \sigma\rangle$.
A property of this bijection which is crucial for us is:
the codimension 1 faces of $\langle \tau^{i-1}\rangle$ are 
$\langle \sigma^i\rangle$, for $\sigma\supseteq\tau$.
Let us pick orientations of all faces $D_T$ of $D$, and extend
them equivariantly to orientations of all cells $\eta$ in $\Sigma_{ghd}$.
Then we orient each dual cell $\ce{\eta}$  dually to
the chosen orientation of $\eta$ (dually with respect to
a chosen orientation of $\Sigma$). These orientations
are of the kind considered in \fullref{sec:4}. 
With these choices we have 
$[\langle\sigma\rangle:\langle\tau\rangle]=
\pm [\sigma:\tau]$, with the sign depending only on 
the dimensions of $\sigma$, $\tau$ (and on $n$, which is fixed in our 
discussion).
    
We define the duality map ${\cal D}\colon L^2_tC^*(\Sigma_{ghd})\to
  L^2_{t^{-1}}C^{n-*}(\Sigma_d)$ by
\begin{equation}
{\cal D}f(\langle\sigma\rangle)=t^{d(\sigma)}f(\sigma).\label{eq:6.2}
\end{equation}
The map ${\cal D}$ is  an isometry of Hilbert spaces.
We will now check that $\delta^{n-i} {\cal D}=\pm {\cal D}\partial^t_i$
(the sign depending only on $i$, $n$):
$$\delta ({\cal D}f)(\langle\tau^{i-1}\rangle)=
\sum_{\sigma^i\supset\tau^{i-1}}[\langle\sigma\rangle:\langle\tau\rangle]
({\cal D}f)(\langle\sigma\rangle)=
\pm\sum_{\sigma^i\supset\tau^{i-1}}[\sigma:\tau]t^{d(\sigma)}f(\sigma)$$      
while
$${\cal D}(\partial^t f)(\langle\tau^{i-1}\rangle)=
t^{d(\tau)}(\partial^t f)(\tau^{i-1})=
t^{d(\tau)}\sum_{\sigma^i\supset\tau^{i-1}}[\sigma:\tau]
t^{d(\sigma)-d(\tau)}f(\sigma)$$
which proves what we wanted.
It follows that ${\cal D}$ intertwines also the adjoint operators;
consequently, it restricts to an isomorphism
${\cal D}\colon L^2_t{\cal H}^*(\Sigma_{ghd})\to L^2_{t^{-1}}
{\cal H}^{n-*}(\Sigma_d)$. 

We still have to check that the Hecke dimensions of these spaces
are the same. 

To this end, let us now consider $L^2_tC^*(\Sigma_{ghd})$
as a subspace of $\oplus_{T\in{\cal F}}L^2_t$ via the embedding
$\Phi_t$ (see \fullref{sec:3}), and 
  $L^2_{t^{-1}}C^{n-*}(\Sigma_d)$ as a subspace of 
$\oplus_{T\in{\cal F}}L^2_{t^{-1}}$ via the embedding
$\Psi_{t^{-1}}$ (see \fullref{sec:4}). We will check that ${\cal D}$
can be regarded as the restriction of the map $D$ (applied
componentwise in $\oplus_{T\in{\cal F}}L^2_t$); it will follow
that $\cal{D}$ preserves dimensions.
Let $f\in L^2(WD_T,\mu_t)$ be a part of a cochain on 
$\Sigma_{ghd}$. Then $$\phi_T(f)=\sqrt{W_T(t)}
\sum_{w\in W^T}f(wD_T)\delta_wp_T(t),$$ where 
$p_T(t)={1\over W_T(t)}\sum_{u\in W_T}
\delta_u$. Since 
{\setlength\arraycolsep{2pt}
\begin{eqnarray}
D(p_T(t))&=&{1\over W_T(t)}\sum_{u\in W_T}(-t)^{d(u)}\delta_u\nonumber\\ 
&=&{1\over W_T((t^{-1})^{-1})}\sum_{u\in W_T}(-t^{-1})^{-d(u)}\delta_u=
h_T(t^{-1}),\nonumber
\end{eqnarray}
}%
we have
\begin{equation}
D(\phi_T(f))=\sum_{w\in W^T}f(wD_T)\sqrt{W_T(t)}(-t)^{d(w)}
\delta_wh_T(t^{-1}).\label{eq:6.3}
\end{equation}
On the other hand,
$({\cal D}f)(w\ce{T})=t^{d(w\ce{T})}f(wD_T)$, and 
\begin{equation}
\psi_T({\cal D}f)=\sum_{w\in W^T}t^{d(w\ce{T})}f(wD_T)(-1)^{d(w)}
\sqrt{W_T(t)}\,\delta_wh_T(t^{-1}).\label{eq:6.4}
\end{equation}
Since for $w\in W^T$ we have $d(w\ce{T})=d(w)$, \eqref{eq:6.3}
and \eqref{eq:6.4} are equal. 
\end{proof}
\begin{remark}
The above proof shows that $\cal{D}$ is an isomorphism of the 
$U_t$--module\break $L^2_t{\cal H}^*(\Sigma_{ghd})$ and the $U_{t^{-1}}$--module $L^2_{t^{-1}}
{\cal H}^{n-*}(\Sigma_d)$, over the algebra isomorphism $D\colon U_t\to U_{t^{-1}}$. 
\end{remark}

\section{Calculation of $b^0_t$}\label{sec:7}

\begin{thm}\label{thm:7.1}
For $t<\rho_W$ we have $b^0_t={1\over W(t)}$; 
for $t\ge\rho_W$ we have $b^0_t=0$.
\end{thm}
\begin{proof} We will use the cell structure $\Sigma_d$.
Vertices of $\Sigma_d$ are located at the centres of chambers $wD$,
thus they are in bijection with $W$. We embed  $L^2_tC^0(\Sigma_d)$
into $L^2_t$ by $(\Psi c)(w)=(-1)^{d(w)}c(w\ce{\emptyset})$. 
This embedding maps all harmonic 0--cochains to constant functions, multiples of ${\bf 1}(w)=1$. 
The square of the
norm of ${\bf 1}$
is $\sum_{w\in W}t^{d(w)}$. It is finite and equal
to $W(t)$ for $t<\rho_W$, and infinite if $t\ge\rho_W$.
The latter means that for $t\ge\rho_W$ we have $L^2_t{\cal H}^0(\Sigma_d)=0$.

To find $b^0_t$ for $t<\rho_W$ we need to identify the projection
of $\delta_1$ on $L^2_t{\cal H}^0(\Sigma_d)$; it is $C{\bf 1}$, where
$$\langle \delta_1-C{\bf 1},{\bf 1}\rangle_t=0.$$
This gives $C=\|{\bf 1}\|^{-2}_t={1\over W(t)}$. In accordance with
the procedure described at the end of \fullref{sec:2}, we find 
$b^0_t=C={1\over W(t)}$.\end{proof}
In view of \fullref{cor:3.4}, the above result makes it 
plausible to suspect that for $t<\rho_W$ we have $b^{>0}_t=0$.
In the next section we prove that this is true for right angled Coxeter groups. 

\section{Mayer--Vietoris sequence}\label{sec:8}

In this section we limit our attention to right angled Coxeter
groups. ``Right angled'' means that whenever two generators $s,s'\in S$ 
are related in the standard presentation, they in fact commute. 
If we join each pair of commuting generators by an edge, we get a graph
with the set of vertices $S$. It is convenient to fill it, 
gluing in a simplex whenever we can
see its 1--skeleton in the graph. The resulting 
simplicial complex is denoted $L$, and the Coxeter group $W_L$.
The Davis chamber $D$ can be identified with the cone $CL'$ over
the first barycentric subdivision of $L$. 
We say that a subcomplex $K\subseteq L$ is full,
if whenever it contains all vertices of a simplex of $L$,
it contains the simplex as well. Full subcomplexes $K$ 
correspond to subsets of $S$ and thus to special subgroups $W_K$
of $W_L$. The Davis complex of $W_K$ is naturally
embedded in $\Sigma_{W_L}$: we first embed
$D_K=CK'$ in $D_L=CL'$, and then extend $W_K$--equivariantly.
We abbreviate $\Sigma_{W_L}$ to $\Sigma_L$.

Let $L=L_1\cup L_2$, where $L_1$, $L_2$ and (consequently) 
$L_0=L_1\cap L_2$ are full subcomplexes of $L$. We embed
$W_{L_i}$ into $W_L$, and $\Sigma_{L_i}$ into $\Sigma_L$;
then $\Sigma_L=W_L\Sigma_{L_1}\cup W_L\Sigma_{L_2}$, 
$W_L\Sigma_{L_1}\cap W_L\Sigma_{L_2}=W_L\Sigma_{L_0}$.
We have a short exact sequence of cochain complexes
$$0\to L^2_tC^*(\Sigma_L)
\to L^2_tC^*(W_L\Sigma_{L_1})\oplus L^2_tC^*(W_L\Sigma_{L_2})\to
L^2_tC^*(W_L\Sigma_{L_0})\to0,$$
from which
we get the long 
Mayer--Vietoris sequence:
\begin{align}
\ldots \to L^2_tH^{i-1}(W_L\Sigma_{L_0})\to& L^2_tH^i(\Sigma_L) \to 
\\
\to L^2_tH^i(W_L\Sigma_{L_1})\oplus& L^2_tH^i(W_L\Sigma_{L_2}) \to
L^2_tH^i(W_L\Sigma_{L_0})\to\ldots\nonumber
\end{align}
Since we work with reduced cohomology,
this sequence is only weakly exact
(the kernels are closures of the images),
see L\"uck \cite[1.22]{L}. Still, if a term is preceded and followed
 by zero terms it has to be zero.
Notice that $W_L\Sigma_{L_i}$ is the disjoint union of 
$w\Sigma_{L_i}$, where $w$ runs through a set of representatives
of $W_{L_i}$--cosets in $W_L$. The $L^2_t$ norm on 
$w\Sigma_{L_i}$ is $t^{d/2}$ times the $L^2_t$ norm on $\Sigma_{L_i}$,
where $d$ is the length of the shortest element of $wW_{L_i}$.
In particular, if $L^2_tH^p(\Sigma_{L_i})=0$, then  
$L^2_tH^p(W_L\Sigma_{L_i})=0$.
\begin{coro}\label{cor:8.1}
Suppose that $b^{>0}_t(\Sigma_{L_i})=0$ for $i=0,1,2$. Then
$b^{>1}_t(\Sigma_L)=0$.
\end{coro}
\begin{thm}\label{thm:8.2}
Let $W$ be a right angled Coxeter group.
For $t<\rho_W$ we have $b^0_t=\chi_t={1\over W(t)}$ and
$b^{>0}_t=0$.
\end{thm}
\begin{proof} Let $W=W_L$. We argue by induction on the number of vertices of $L$.

(1)\qua If $L$ is a simplex, then $\Sigma_{L,d}$ is a cube; its $L^2_t$ 
cohomology coincides with the usual cohomology and is concentrated 
in dimension 0. 

(2)\qua If $L$ is not a simplex, we can find two vertices $a,b\in L$
not connected by an edge; we put $L_1=\bigcup\{\sigma\mid
a\not\in\sigma\}$, $L_2=\bigcup\{\sigma\mid
b\not\in\sigma\}$ and $L_0=L_1\cap L_2$. 
These have fewer vertices than $L$, and so
$L^2_tH^{>0}(\Sigma_{L_i})=0$ for $t<\rho(W_{L_i})$ ($i=0,1,2$). 
Since $L_i\subset L$, we have $\rho(W_{L_i})\ge\rho(W_L)$. Therefore
we have $L^2_tH^{>0}(\Sigma_{L_i})=0$ for $t<\rho(W_L)$.
It follows from \fullref{cor:8.1} that $L^2_tH^{>1}(\Sigma_L)=0$ (still for 
$t<\rho(W_L)$), while from \fullref{cor:3.4} and \fullref{thm:7.1} we conclude
that
$$b^0_t(\Sigma_L)=\chi_t(\Sigma_L)=b^0_t(\Sigma_L)-
b^1_t(\Sigma_L).$$
Thus $b^1_t(\Sigma_L)=0$.\end{proof}
\begin{coro}\label{cor:8.3}
Assume that $L$ is a generalised homology $(n-1)$--sphere
(ie, $(D,\partial D)$ is a generalised homology $n$--disc); then for
$t<{1\over \rho(W_L)}$ we have $b^n_t=0$,
while for $t>{1\over \rho(W_L)}$ the $L^2_t$--cohomology
is concentrated in dimension $n$ and $b^n_t=(-1)^n\chi_t={(-1)^n\over W_L(t)}$.
\end{coro}
\begin{proof} This follows from Theorems \ref{thm:8.2} and \ref{thm:7.1} via Poincar\' e duality
(\fullref{thm:6.1}). 
\end{proof}
\begin{propo}\label{prop:8.4}
Let $K\subset L$ be a full subcomplex. The dimension
of the $U_t(W_L)$--module $L^2_tH^q(W_L\Sigma_K)$ is the same
as the dimension of the $U_t(W_K)$--module\break $L^2_tH^q(\Sigma_K)$
(ie, it is equal to $b^q_t(\Sigma_K)$).
\end{propo}
\begin{proof} 
A harmonic $q$--cochain on $W_L\Sigma_K=\bigcup\{w\Sigma_K\mid w\in W_L\}$
is the same thing as a collection of harmonic $q$--cochains on $w\Sigma_K$.
In order to calculate dimensions, we embed everything in $V=\oplus_{\sigma\subset D_L}L^2_t(W_L)$.
Let  ${\bf 1}_{\sigma}\in V$ have $\delta_1$ as its coordinate with index $\sigma$,
and 0 on all other coordinates.   
As we project ${\bf 1}_{\sigma^q}$ on 
$L^2_t{\cal H}^q(W_L\Sigma_K)$,
we get in fact a harmonic cochain supported on $\Sigma_K$---harmonic 
cochains supported on other components of $W_L\Sigma_K$
are orthogonal to ${\bf 1}_{\sigma^q}$, so also to its projection.
We can as well project ${\bf 1}_{\sigma^q}$ on $L^2_t{\cal H}^q(\Sigma_K)$
inside $\oplus L^2_t(W_K)$, so that the projection matrices are the same
(apart for the case $\sigma\not\subset K$, which
gives 0 in the first case and does not appear in the second),
and traces coincide.\end{proof}

\section{Chain homotopy contraction}\label{sec:9}

In this section we will describe a simplicial version of the 
geodesic contraction of $\Sigma$ with respect to the Moussong metric.
We will consider the chain complex $C_*(\Sigma_{st})$ equipped with the 
boundary operator $\partial$ given by \eqref{eq:5.7}.
Henceforth we write $\Sigma$ for $\Sigma_{st}$,
and we denote by  $b$  the barycentre of the basic chamber $D$. 
Recall that $\Sigma$ can be equipped with a $W\!$--invariant, $CAT(0)$ metric $d_M$, the \metric (Moussong \cite{M}).
From now on, all balls, geodesics etc. will be considered with respect to $d_M$ (unless explicitly
stated otherwise).
Besides $CAT(0)$, the following property of the \metric will be useful for us:
for every $R>0$ there exists a constant $N(R)$ such that 
any ball of radius $R$ in $\Sigma$ intersects at most $N(R)$ chambers.
\begin{thm}\label{thm:9.1}
There exists a linear map $H\colon C_*(\Sigma)\to C_{*+1}(\Sigma)$, and constants $C$, $R$,
with the following properties:
\begin{itemize}
\item[\rm(a)] if $v\in\Sigma^{(0)}$, then $\partial H(v)=v-b$;
\item[\rm(b)] if $\sigma$ is a simplex of positive dimension, then $\partial H(\sigma)=\sigma
-H(\partial \sigma)$;
\item[\rm(c)] for every simplex $\sigma$, $\|H(\sigma)\|_{L^\infty}<C$;
\item[\rm(d)] if $\gamma$ is a geodesic from a vertex
of a simplex $\sigma$ to $b$, then $\supp(H(\sigma))\subseteq B_R({\rm image}(\gamma))$.
\end{itemize}
\end{thm}
\begin{proof} 
We will construct, for all integers $i\ge0$, linear maps
$h_i\colon C_*(\Sigma)\to C_*(\Sigma)$, $H_i\colon C_*(\Sigma)\to C_{*+1}(\Sigma)$ such that:
\begin{itemize}
\item[(1)] $h_0={\rm id}$;
\item[(2)] $\partial h_i=h_i\partial$;
\item[(3)] $\partial H_i=h_i-H_i\partial-h_{i+1}$;
\item[(4)] $\exists C_k,\forall\sigma\in\Sigma^{(k)},\forall i\ge0,
\|H_i(\sigma)\|_{L^\infty}<C_k$ 
and 
$\|h_i(\sigma)\|_{L^\infty}<C_k$; 
\item[(5)] $\exists R_k,\forall\sigma\in\Sigma^{(k)},\forall i\ge0$, if $\gamma$ is a geodesic from a vertex of $\sigma$
to $b$, then $\supp(h_i(\sigma))$, $\supp(H_{i-1}(\sigma))$ (if $i>0$) 
and $\supp(H_i(\sigma))$ are
contained in the ball $B_{R_k}(\gamma(i))$ (or in $B_{R_k}(b)$, if $i>{\rm length}(\gamma)$);  
\item[(6)] if $i\ge{\rm diam}(\sigma\cup\{b\})$, then $h_i(\sigma)=0$ (unless $\dim\sigma=0$, in which case
$h_i(\sigma)=b$) and $H_i(\sigma)=0$.
\end{itemize}
 
The construction will be by induction on the chain 
degree $k$.
Throughout this proof, we will say that a family of chains is uniformly
bounded if they have uniformly bounded support
diameters and $L^\infty$ norms. Let $A$ be the length of the longest 
edge in $\Sigma$.
\eject

(1)\qua $k=0$\par
Let $v\in\Sigma^{(0)}$, let $\gamma_v\colon [0,l]\to\Sigma$ be a geodesic such that
$\gamma_v(0)=v$, $\gamma_v(l)=b$. We put $h_0(v)=v$, $h_i(v)=b$ if $i\ge l$, and we choose
a vertex within distance $A$ from $\gamma_v(i)$ and declare it to be $h_i(v)$ in the remaining cases.  
We have $d(h_i(v), h_{i+1}(v))\le 1+2A$. Now, up to the action of
$W$, there are only finitely many pairs of vertices $(y,z)$ satisfying $d(y,z)<1+2A$.
In every $W\!$--orbit of such pairs we choose a pair $(y,z)$ 
and we fix a 1--chain $H(y,z)$, $\partial H(y,z)=y-z$;
we then extend $H$ to the $W\!$--orbit of $(y,z)$ using the $W\!$--action
(making choices if stabilisers are non-trivial).
In the case $y=z$ we choose $H(y,y)=0$.
Notice that the chosen 1--chains $H$ are uniformly bounded.
Finally, we put $H_i(v)=H(h_i(v),h_{i+1}(v))$.  

(2)\qua $k\to(k+1)$\par
Let $\sigma\in\Sigma^{(k+1)}$. Then, due to (2), $\partial h_i(\partial\sigma)=h_i(\partial\partial\sigma)=0$.
Thus, $h_i(\partial\sigma)$ is a cycle. Moreover, we claim that as we vary $\sigma$, the cycles $h_i(\partial\sigma)$
are uniformly bounded.
In fact, as a consequence of (5), every simplex in the support of $h_i(\partial\sigma)$
is within $R_k$ of one of the points $\gamma_v(i)$, where
$v$ runs through the vertices of $\sigma$, and, by $CAT(0)$ comparison, 
the $k+2$ points $\gamma_v(i)$ are within 
$2A$ of each other. Whence uniform boundedness of supports of $h_i(\partial\sigma)$.
Uniform boundedness of $L^\infty$ norms follows from (4).  
Up to the $W\!$--action on $C_k(\Sigma)$, there are only finitely many possible
values of $h_i(\partial\sigma)$. 
As in step 1, we fix $(k+1)$--chains $h_i(\sigma)$,
$\partial h_i(\sigma)=h_i(\partial\sigma)$, so that they are uniformly bounded
(and are 0 whenever $h_i(\partial\sigma)=0$).

To define $H_i(\sigma)$, we consider the chain $h_i(\sigma)-H_i(\partial\sigma)-h_{i+1}(\sigma)$.
It is a cycle:
$$\partial(h_i(\sigma)-H_i(\partial\sigma)-h_{i+1}(\sigma))=
\partial h_i(\sigma)-\partial H_i(\partial\sigma)-\partial h_{i+1}(\sigma)$$
$$= 
h_i(\partial\sigma)-\bigl(h_i(\partial\sigma)-H_i(\partial\partial\sigma)-h_{i+1}(\partial\sigma)\bigr)
-h_{i+1}(\partial\sigma)=0.$$
Again, all such chains (as we vary $\sigma$) are uniformly bounded, 
and we can choose $H_i(\sigma)$, satisfying $\partial H_i(\sigma)=h_i(\sigma)-H_i(\partial\sigma)-h_{i+1}(\sigma)$,
in a uniformly bounded way. As before, we put $H_i(\sigma)=0$ whenever we have to chose it 
so that it has boundary 0 (so as to satisfy (6)). 

Now that we have a family of maps satisfying (1)--(6), we put $H(\sigma)=\sum_{i\ge0}H_i(\sigma)$.
The sum is always finite because of (6). The conditions (a)--(d) are easy to check:
(a) and (b) follow from (1), (3) and (6);
(c) follows from (4) and (5): since the supports of $H_i(\sigma)$ are uniformly bounded and
``move along'' a geodesic $\gamma$ with constant speed as $i$ grows, 
only a uniformly finite number of $H_i(\sigma)$
contribute to a coefficient of a fixed simplex $\tau$ in the chain $H(\sigma)$; moreover,
because of (4), each contribution is smaller than $C_{\dim{\sigma}}$;
(d) is a consequence of (5).\end{proof}

\section{Vanishing below $\rho$}\label{sec:10}

Let $H$ be a map as in \fullref{thm:9.1}.
\begin{thm}\label{thm:10.1}
Suppose that $t>{1\over\rho_W}$. Then the map $H$ extends to 
a bounded operator $H\colon L^2_tC_*(\Sigma)\to L^2_tC_{*+1}(\Sigma)$.
\end{thm}
\begin{proof} 
Unspecified summations will be over $\Sigma^{(k)}$. $N_k$ will denote the number of
$k$--simplices in a chamber.  

Let $a=\sum a_\sigma\sigma\in L^2_tC_k(\Sigma)$. We know that for every simplex $\sigma$,
$\|H(\sigma)\|_{L^\infty}<C$. Also
{\setlength\arraycolsep{2pt}
\begin{eqnarray}
\sum |a_\sigma|&=&\sum |a_\sigma|t^{d(\sigma)/2}t^{-d(\sigma)/2}\le
\left(\sum|a_\sigma|^2t^{d(\sigma)}\right)^{1/2}\left(\sum t^{-d(\sigma)}\right)^{1/2}\nonumber\\
&\le&   
\|a\|_t\left(N_kW(t^{-1})\right)^{1/2}<+\infty,\nonumber
\end{eqnarray}
}%
so that $\sum a_\sigma H(\sigma)$ is pointwise convergent to a chain 
$H(a)\in L^\infty C_{k+1}(\Sigma)$. We want to estimate $\|H(a)\|_t$.
Let us write $\tau\prec\sigma$ if $\tau$ appears with non-zero coefficient in
$H(\sigma)$. We have $|H(a)_\tau|\le\sum_{\sigma\mid \tau\prec\sigma}C|a_\sigma|$
, so that
{\setlength\arraycolsep{2pt}
\begin{eqnarray}
\sum|&&\hbox{\hskip-.47cm}H(a)_\tau|^2t^{d(\tau)}\le  
C^2\sum_\tau\biggl(\sum_{\sigma\mid \tau\prec\sigma}|a_\sigma|\biggr)^2t^{d(\tau)}\nonumber\\
&\le& C^2\sum_\tau \biggl(\sum_{\sigma\mid \tau\prec\sigma}|a_\sigma| t^{d(\sigma)/2}     
t^{-\alpha\left({d(\sigma)-d(\tau)\over2}\right)}
t^{-\beta\left({d(\sigma)-d(\tau)\over2}\right)}\biggr)^2\label{eq:10.1}\\
&\le&
C^2\sum_\tau\biggl(\sum_{\sigma\mid \tau\prec\sigma}|a_\sigma|^2 t^{d(\sigma)}     
(t^{-\alpha})^{d(\sigma)-d(\tau)}\biggr)
\biggl(\sum_{\sigma\mid \tau\prec\sigma}
(t^{-\beta})^{d(\sigma)-d(\tau)}\biggr).\nonumber
\end{eqnarray}
}%
Here $\alpha$, $\beta$ are positive numbers chosen so that $\alpha+\beta=1$, $t^{-\beta}<\rho_W$.
\begin{claim}
There exists a constant $C'$, independent of $\tau$, such that 
$$\sum_{\sigma\mid \tau\prec\sigma}
(t^{-\beta})^{d(\sigma)-d(\tau)}\le C'W(t^{-\beta}).$$
\end{claim}
\begin{proof} Recall that $A$ is the length of the longest edge in $\Sigma$,
and $N(r)$ is the maximal number of chambers intersecting a ball of radius $r$.
The claim follows from two observations.

(1)\qua For $w_0\in W$ let $E(w_0)=\{w\in W\mid d(w)=d(w_0)+d(w_0^{-1}w)\}$. In more geometric terms,
$E(w_0)$ is the set of all $w$ such that some gallery connecting $D$ and $wD$ passes through 
$w_0D$. We have
$$\sum_{w\in E(w_0)}(t^{-\beta})^{d(w)-d(w_0)}=
 \sum_{w\in E(w_0)}(t^{-\beta})^{d(w_0^{-1}w)}\le
\sum_{w\in W}(t^{-\beta})^{d(w)}=W(t^{-\beta}).$$

(2)\qua If $\tau\prec\sigma$, then $\tau$ is at distance at most $R$ from 
a geodesic $\gamma$ joining (a vertex of) $\sigma$ and $b$. 
Let us consider the union $U$ of all galleries joining $D$ and a fixed chamber $D'$
containing $\sigma$. Then $U$ is the intersection of all half-spaces
containing $D$ and $D'$ (see Ronan \cite{Ron}). 
Since half-spaces are geodesically convex in $d_M$, we have
$\gamma\subseteq U$. 
Consequently, every point of $\gamma$ lies in a 
gallery joining $D'$ and $D$. 
Therefore, if we put $B(\tau)=\{w_0\mid w_0D\cap B_R(\tau)\ne\emptyset\}$, then we have
$\{\sigma\mid \tau\prec\sigma\}\subseteq\bigcup_{w_0\in B(\tau)}E(w_0)D$.

Putting these together,
{\setlength\arraycolsep{2pt}
\begin{eqnarray}
\sum_{\sigma\mid \tau\prec\sigma}
(t^{-\beta})^{d(\sigma)-d(\tau)}
&\le& 
\sum_{w_0\in B(\tau)}t^{-\beta(d(w_0)-d(\tau))}
\sum_{w\in E(w_0)}N_k(t^{-\beta})^{d(w)-d(w_0)}\nonumber\\
&\le&
\sum_{w_0\in B(\tau)}t^{-\beta(d(w_0)-d(\tau))}
N_kW(t^{-\beta}).\nonumber
\end{eqnarray}
}%
Notice that $|d(w_0)-d(\tau)|$ does not exceed the gallery
distance from $w_0D$ to some chamber containing $\tau$, and is therefore
uniformly bounded. Also, the cardinality of $B(\tau)$ is
bounded by $N(R+A)$. The claim is proved.
\end{proof}
Using the claim, we can continue the estimate \eqref{eq:10.1}:
$$\|H(a)\|_t^2\le C^2C'W(t^{-\beta})
\sum_\tau\biggl(\sum_{\sigma\mid \tau\prec\sigma}|a_\sigma|^2 t^{d(\sigma)}
(t^{-\alpha})^{d(\sigma)-d(\tau)}\biggr).$$ 
Now 
$$\sum_\tau\biggl(\sum_{\sigma\mid \tau\prec\sigma}|a_\sigma|^2 t^{d(\sigma)}     
(t^{-\alpha})^{d(\sigma)-d(\tau)}\biggr)=\sum_\sigma\biggl(|a_\sigma|^2 t^{d(\sigma)}     
\sum_{\tau\mid \tau\prec\sigma}
(t^{-\alpha})^{d(\sigma)-d(\tau)}\biggr),$$
so that the following lemma is all we need:
\begin{lem}\label{lem:10.2}
There exists a 
constant $K$ independent of $\sigma$ such that 
$$\sum_{\tau\mid \tau\prec\sigma}
(t^{-\alpha})^{d(\sigma)-d(\tau)}<K.$$
\end{lem}
\proof Since $W$ acts on $(\Sigma,d_M)$ isometrically, cocompactly and properly discontinuously,
the word metric $d$ on $W$ is quasi-isometric to the metric $d_M$ restricted to 
$W\simeq Wb\hookrightarrow\Sigma$. This implies that there exist constants $M$, $m$, $L$ such that
for any two points $y,z\in\Sigma$ and any chambers $D_y\ni y$, $D_z\ni z$, we have
\begin{equation}
Md_M(y,z)+L\ge d(D_y,D_z)\ge md_M(y,z)-L,\label{eq:10.2}
\end{equation}
where we put $d(wD,uD)=d(w,u)=d(w^{-1}u)$.

Let $v$ be a vertex of $\sigma$, and let $\gamma\colon [0,l]\to\Sigma$ be a geodesic,
$\gamma(0)=v$, $\gamma(l)=b$. To each $\tau\prec\sigma$ we can assign one of the points $\gamma(i)$ 
($0\le i\le\lfloor l\rfloor$)  
in such a way that $d_M(\tau,\gamma(i))<R+1$. The number
of simplices to which we assign a given $\gamma(i)$ does not exceed 
$N(R+1)N_k$.
Suppose that $\gamma(i)$ is assigned to $\tau$. Let $D_\tau$ 
(resp.\ $D_\sigma$) be the chamber containing $\tau$ 
(resp.\ $\sigma$) 
such that
$d(\tau)=d(D,D_\tau)$ (resp.\ $d(\sigma)=d(D,D_\sigma)$). 
Let $D_i$ be a chamber containing $\gamma(i)$. 
We choose $D_i$ so that some gallery from $D$ to $D_\sigma$ passes through $D_i$
(see part 2 of the proof of the claim above).
Using \eqref{eq:10.2}, we get 
{\setlength\arraycolsep{2pt}
\begin{eqnarray}
d(\sigma)-d(\tau)&=&d(D,D_\sigma)-d(D,D_\tau)\nonumber\\ 
&\ge&
d(D,D_i)+d(D_i,D_\sigma)-(d(D,D_i)+d(D_i,D_\tau))\nonumber\\ 
&\ge&
md_M(\gamma(i),v)-L-(Md_M(\tau,\gamma(i))+L)\nonumber\\  
&\ge& 
mi-(M(R+1)+2L)=mi-P,\nonumber
\end{eqnarray}
}%
where $P=M(R+1)+2L$. 
Remember that $t^{-1}$ and, whence, $t^{-\alpha}$ are less than~1.
Therefore
{\setlength\arraycolsep{2pt}
\begin{eqnarray}
\sum_{\tau\mid \tau\prec\sigma}
(t^{-\alpha})^{d(\sigma)-d(\tau)}&\le&
\sum_{i=0}^{\lfloor l\rfloor} N(R+1)N_k(t^{-\alpha})^{mi-P}\nonumber\\
&=& 
N(R+1)N_kt^{\alpha P}\sum_{i=0}^{\lfloor l\rfloor}(t^{-\alpha m})^i\nonumber\\
&\le&
N(R+1)N_kt^{\alpha P}{1\over 1-t^{-\alpha m}}.\nonumber
\end{eqnarray}
}%
This completes the proof of \fullref{lem:10.2} and of \fullref{thm:10.1}.
\end{proof}
\begin{thm}\label{thm:10.3}
Let $W$ be a Coxeter group.
For $t<\rho_W$ we have $b^0_t(W)=\chi_t(W)={1\over W(t)}$ and
$b^{>0}_t(W)=0$.
\end{thm}
\begin{proof} Theorems \ref{thm:9.1} and \ref{thm:10.1} imply that in the range $t>{1\over\rho_W}$
we have
$$H_{>0}(L^2_tC_*(\Sigma),\partial)=0.$$ Indeed, if 
$c\in L^2_tC_k(\Sigma)$, $\partial c=0$, then
$c=\partial H(c)+H(\partial c)=\partial H(c)$, so that $[c]=0$.
It follows that the isomorphic complex
$(L^2_{t^{-1}}C_*(\Sigma),\partial^{t^{-1}})$ also has vanishing
homology in degrees $>0$ (if $t^{-1}<\rho_W$). 
Thus, its homology is concentrated in dimension 0,
and the zeroth Betti number is equal to the Euler characteristic.
\end{proof}
\begin{coro}\label{cor:10.4}
Assume that $(D,\partial D)$ is a generalised homology $n$--disc; then for
$t<{1\over \rho_W}$ we have $b^n_t=0$,
while for $t>{1\over \rho_W}$ the $L^2_t$ cohomology
is concentrated in dimension $n$ and $b^n_t=(-1)^n\chi_t={(-1)^n\over W(t)}$.
\end{coro}
\begin{proof} This follows from Theorems \ref{thm:10.3} and \ref{thm:7.1} using Poincar\' e duality
(\fullref{thm:6.1}).
\end{proof}

\bibliographystyle{gtart}
\bibliography{link}

\end{document}